\documentclass[reqno]{amsart} 

\begin{document} 
\newtheorem{theorem}{Theorem}[section]
\newtheorem{lemma}{Lemma}[section]
\newtheorem{corollary}[theorem]{Corollary}
\newtheorem{example}[theorem]{Example}
\newtheorem{remark}{Remark}{}
\title[\hfil Remarks on an 1828 theorem of Clausen]
{Remarks on an 1828 theorem of Clausen} 

\author[Angelo B. Mingarelli \hfilneg]
{Angelo B. Mingarelli$^{(1)}$}  

\address{School of Mathematics and Statistics\\ 
Carleton University, Ottawa, Ontario, Canada, K1S\, 5B6.}
\email[A. B. Mingarelli]{amingare@math.carleton.ca, amingarelli@dma.ulpgc.es}
\address{Departamento de Matem\'aticas, Universidad de Las Palmas de Gran Canaria, 
Campus de Tafira 
Baja, \\ 35017 Las Palmas de Gran Canaria, Spain.}

\date{}
\thanks{Preprint [August, 2006]. Revised June 18, 2008.}
\thanks{The author is partially supported by a NSERC Canada Research Grant. Gratitude is also expressed to the Department of Mathematics of the University of Las Palmas in Gran Canaria for its hospitality during the author's term there as a Visiting Professor, 2006-2008}

\begin{abstract}
We list some explicit calculations related to a theorem of Clausen originally published in 1828, more commonly known as the result that describes the linear third order differential equation satisfied by the squares and the product of any two solutions of a linear second order differential equation in the real domain. The case of the cube of a solution dates to Appell, 1880. Although not commonly known and perhaps even new to some extent we show that the fourth and fifth powers of such solutions actually satisfy a linear differential equation of order five and six respectively provided the coefficients are sufficiently smooth. Indeed, it is the case that given any solution of a linear second order equation its m-th power satisfies an effectively computable linear differential equation of order m+1. 
\end{abstract}

\maketitle

\section{The differential equations}

Consider the real linear second order differential equation in an unspecified interval $I$, of the form
\begin{equation}
\label{eq1}
y^{\prime\prime}=p(x)y^{\prime} + q(x)y, \quad \quad x\in I,
\end{equation}
where $p, q$ are sufficiently smooth functions to allow the operations below to proceed. For example, if $f$ is a solution of \eqref{eq1} and we are seeking the differential equation satisfied by the function $y(x)=f(x)^m$, $m \geq 3$, then $p, q \in C^{m-1}(I)$. In what follows two linearly independent solutions of \eqref{eq1} will always be denoted by $f,g$. The case $m=2$ is well-known and admits a myriad of applications: The result can be stated simply as saying that the quantities $f^2, fg, g^2$ form a basis for a related linear third order differential equation, that is, cf., [\cite{WW}, p.298]
\begin{equation}
\label{eq2}
y^{\prime\prime\prime}=3p\, y^{\prime\prime} + (4q - 2p^2 + p^{\prime})\, y^{\prime} + 2(q^{\prime} - 2pq)\, y, \quad \quad x\in I,
\end{equation}
\noindent Its early appearance can be found in the books by Ince \cite{ince} and Whittaker and Watson \cite{WW}, although the latter credit Appell \cite{app}, the case $m=2$ appears to be due to Clausen \cite{cl}, see Gasper and Rahman [\cite{gr}, p.xv]. In the developments stated below all proofs are lengthy but straightforward and can be ascertained by direct verification. The third power analog of this result is as follows: Whenever $f^{\prime}(x)\neq 0$ and $g^{\prime}(x)\neq 0$, we have the following results:

\begin{theorem}The functions $f^3, f^2g, fg^2, g^3$ form a basis of the linear fourth order differential equation
\label{th1}
\begin{eqnarray}
\label{eq3}
&&y^{(iv)}- 6\,p \,y^{\prime\prime\prime} + (11p^2-4p^{\prime}-10q)\,y^{\prime\prime}+ (7pp^{\prime}+30pq-6p^3-p^{\prime\prime}-10q^{\prime})\,y^{\prime}\nonumber \\
&& + (9q^2+15pq^{\prime}+6p^{\prime}q-18p^2q-3q^{\prime\prime})\, y = 0.
\end{eqnarray}
\end{theorem}

\noindent The case $m=3$ above is due to Appell~[\cite{app},p.213] but is not commonly known; none of the classic authors appear to mention it or any other generalizations. We have not been able to find any references for the next two cases, nevertheless they are interesting enough to be listed since the study of orthogonal polynomials has reached equations of degree at least six these days.

\vskip0.15in

\begin{theorem}The functions $f^4, f^3g, f^2g^2, fg^3,g^4$ form a basis of the linear fifth order differential equation
\label{th1}
\begin{eqnarray}
\label{eq4}
&&y^{(v)} + c_4 \,y^{(iv)} + c_3\,y^{\prime\prime\prime}+ c_2\,y^{\prime\prime}+ c_1\, y^{\prime}+c_0\,y = 0.
\end{eqnarray}
\end{theorem}
\noindent where \\
\noindent $c_4=-10p,$\\
\noindent $c_3 = 35p^2 - 10p^{\prime} -20q,$\\
\noindent $c_2 = 4pp^{\prime} - 5 p^{\prime\prime} -30q^{\prime} +120pq - 50p^3,$\\
\noindent $c_1 = 64q^2+24p^4+56p^{\prime}q+11pp^{\prime\prime}-46p^2p^{\prime}-208p^2q+128pq^{\prime}-p^{\prime\prime\prime}+7{p^{\prime}}^2-18q^{\prime\prime},$\\
\noindent $c_0 = 36pq^{\prime\prime}-128pq^2-4q^{\prime\prime\prime}+28p^{\prime}q^{\prime}-104p^2q^{\prime}+64qq^{\prime}-80pp^{\prime}q+8p^{\prime\prime}q+96p^3q.$

\vskip0.15in

\begin{theorem}The functions $f^5, f^4g, f^3g^2, f^2g^3,fg^4,g^5$ form a basis of the linear sixth order differential equation
\label{th1}
\begin{eqnarray}
\label{eq4}
&&y^{(vi)} + c_5y^{(v)}+ c_4 \,y^{(iv)} + c_3\,y^{\prime\prime\prime}+ c_2\,y^{\prime\prime}+ c_1\, y^{\prime}+c_0\,y = 0.
\end{eqnarray}
\end{theorem}
\noindent where \\
\noindent $c_5 = -15p,$\\
\noindent $c_4=  85p^{2}-35q -20p^{\prime},$\\
\noindent $c_3 = 350qp - 225p^3 -70q^{\prime} - 15p^{\prime\prime} +165pp^{\prime},$\\
\noindent $c_2 = 274p^4 - 63q^{\prime\prime} - 1183p^2q + 52 {p^{\prime}}^2 +525pq^{\prime}+81pp^{\prime\prime}-421p^2p^{\prime}+259q^2+266qp^{\prime}-6p^{\prime\prime\prime},$\\
\noindent $c_1 = 315pq^{\prime\prime} -120p^5-28q^{\prime\prime\prime}-p^{(iv)} -1071pp^{\prime}q +326p^3p^{\prime}-1295pq^2+1540p^3q-101p^2p^{\prime\prime}+ 16pp^{\prime\prime\prime} + 25p^{\prime}p^{\prime\prime}+91qp^{\prime\prime}+266p^{\prime}q^{\prime}-1183p^2q^{\prime} - 127p{p^{\prime}}^2 + 518qq^{\prime},$\\
\noindent $c_0 = 1450q^2p^2 - 1295qq^{\prime}p - 310q^2p^{\prime}-5q^{(iv)}+860p^2p^{\prime}q- 535pp^{\prime}q^{\prime}-150pp^{\prime\prime}q +770p^3q^{\prime} - 600p^4q - 100{p^{\prime}}^2q -355p^2q^{\prime\prime}+45p^{\prime\prime}q^{\prime}+80p^{\prime}q^{\prime\prime}+10p^{\prime\prime\prime}q + 155qq^{\prime\prime} + 70pq^{\prime\prime\prime}+130{q^{\prime}}^2-225q^3.$ \\

\noindent {\bf Remark} An (m+1)-th order linear ordinary differential equation with the basis elements $f^{m},f^{m-1}g,\ldots,f^ig^j,\ldots, g^{m}$ where $i+j=m$, can be derived as well. Thus, higher order equations of the type presented above can be integrated easily.


\begin{thebibliography}{+99}
\bibitem{app}P. Appell, {\it Sur la transformation des \'{e}quations diff\'{e}rentielles lin\'{e}aires}, in Comptes Rendus Acad. Sci. Paris, $2^e$ Semestre {\bf 91} (4) (1880), 211-214
\bibitem{cl} T. Clausen, {\it \"{U}ber die F\"{a}lle \ldots}, J. Reine Ang. Math. [Crelle's Journal] {\bf 3}, (1828), 89-91
\bibitem{gr} G. Gasper and M. Rahman, {\it Basic Hypergeometric Series}, Second Edition, Cambridge University Press, Cambridge, UK. 2004
\bibitem{ince} E.L. Ince,\,\,{\it Ordinary Differential Equations},
Dover Publications, New York, 1956
\bibitem{WW} E.T.Whittaker and G.N.Watson, {\it A Course of Modern Analysis}, Fourth Edition, Cambridge University Press, Cambridge, UK., 1927
\end{thebibliography}
\end{document}